\documentclass[11pt]{amsart}
\usepackage{amsmath,amssymb}
\usepackage{curves,eepicemu, eepic}
\usepackage{amsmath,amssymb}

\newbox\TempBox \newbox\TempBoxA
\def\UnderWiggly#1{
  \ifmmode\setbox\TempBox=\hbox{$#1$}\else\setbox\TempBox=\hbox{#1}\fi
  \setbox\TempBoxA=\hbox to \wd\TempBox{\hss\char'176\hss}
  \rlap{\copy\TempBox}\smash{\lower12pt\hbox{\copy\TempBoxA}}}

\newtheorem{theorem}{Theorem}[section]
\newtheorem{proposition}[theorem]{Proposition}
\newtheorem{corollary}[theorem]{Corollary}
\newtheorem{conjecture}[theorem]{Conjecture}
\newtheorem{lemma}[theorem]{Lemma}
\newtheorem{claim}[theorem]{Claim}

\newtheorem{remark}[theorem]{Remark}
\newtheorem{defi}[theorem]{Definition}

\newcommand{\bt}{\begin{theorem}}
\newcommand{\et}{\end{theorem}}
\newcommand{\bl}{\begin{lemma}}
\newcommand{\el}{\end{lemma}}
\newcommand{\bcl}{\begin{claim}}
\newcommand{\ecl}{\end{claim}}
\newcommand{\eq}[1]{\begin{equation} #1 \end{equation}}
\newcommand{\bp}{\begin{proposition}}
\newcommand{\ep}{\end{proposition}}
\newcommand{\bcor}{\begin{corollary}}
\newcommand{\ecor}{\end{corollary}}
\newcommand{\br}{\begin{remark}\rm}
\newcommand{\er}{\end{remark}}
\newcommand{\bcon}{\begin{conjecture}}
\newcommand{\econ}{\end{conjecture}}
\def\push{\hspace{1in}}

\def\la{\langle}
\def\ra{\rangle}
\def\dmas{D([0,\infty), \mfs)}
\def\demas{D([\epsilon,\infty), \mfs)}
\def\qed{\hfill$\Box$}

\newcommand{\ben}{\begin{equation}}
\newcommand{\een}{\end{equation}}
\newcommand{\be}{\begin{enumerate}}
\newcommand{\ee}{\end{enumerate}}
\newcommand{\bei}{\begin{itemize}}
\newcommand{\eei}{\end{itemize}}
\newcommand{\bea}{\begin{eqnarray}}
\newcommand{\eea}{\end{eqnarray}}
\newcommand{\beas}{\begin{eqnarray*}}
\newcommand{\eeas}{\end{eqnarray*}}
\newcommand{\ba}{\begin{array}}
\newcommand{\ea}{\end{array}}
\newcommand{\bc}{\be\begin{array}{r@{\,}c@{\,}l}}
\newcommand{\ec}{\end{array}\ee}

\def\mfs{M({\mathbb R}_+ \times {\mathbb R})}
\def\mfas{M_a({\mathbb R}_+ \times {\mathbb R})}
\newcommand{\al}{\alpha}

\newcommand{\eps}{\varepsilon}

\newcommand{\f}{\frac}

\def\cid{\stackrel{d}{\longrightarrow}}

\def\eid{\stackrel{d}{=}}

\def\nmr{\nonumber}

\newcommand{\Bi}{{\mathcal B}}
\newcommand{\Ci}{{\mathcal C}}

\newcommand{\Ei}{{\mathcal E}}

\newcommand{\Gi}{{\mathcal G}}

\newcommand{\Li}{{\mathcal L}}

\newcommand{\Ri}{{\mathcal R}}

\newcommand{\Ti}{{\mathcal T}}

\newcommand{\Yi}{{\mathcal Y}}

\newcommand{\R}{{\mathbb R}}
\newcommand{\N}{{\mathbb N}}

\newcommand{\querry}[1]{}

\makeatletter\@addtoreset{equation}{section}
\makeatother

\begin{document}

\author{Krishna B. Athreya}
\address{Krishna B. Athreya \\ Department of Mathematics and Statistics\\ Iowa State University\\ Ames, Iowa, 50011, U.S.A.}
\email{kba@iastate.edu}

\author{Siva R. Athreya}
\address{Siva R. Athreya \\ 8th Mile Mysore Road\\ Indian Statistical Institute
  \\Bangalore 560059, India.}
\email{athreya@isibang.ac.in}

\author{Srikanth K. Iyer}
\address{Srikanth K. Iyer\\ Department of Mathematics\\ Indian Institute of Science\\ Bangalore  560012, India.}
\email{skiyer@math.iisc.ernet.in}


\title[Age Dependent Branching Markov Processes]{Critical Age Dependent Branching Markov Processes and their Scaling Limits}

\begin{abstract}
This paper studies: (i) the long time behaviour of the empirical
distribution of age and normalised position of an age dependent
critical branching Markov process conditioned on non-extinction;
and (ii) the super-process limit of a sequence of age dependent
critical branching Brownian motions.
\end{abstract}

\maketitle

\vspace{-.7cm}

\section{Introduction}

Consider an age dependent branching Markov process where i) each
particle lives for a random length of time and during its lifetime
moves according to a Markov process and ii) upon its death it
gives rise to a random number of offspring. We assume that the
system is critical, i.e. the mean of the offspring distribution
is one.

We study three aspects of such a system. First, at time $t,$
conditioned on non-extinction (as such systems die out w.p. $1$)
we consider a randomly chosen individual from the population. We
show that asymptotically (as $t \rightarrow \infty$), the joint
distribution of the position (appropriately scaled) and age
(unscaled) of the randomly chosen individual decouples (See
Theorem \ref{oneparticle}). Second, it is shown that conditioned
on non-extinction at time $t,$ the empirical distribution of the
age and the normalised position of the population converges as $t
\rightarrow \infty$ in law to a random measure characterised by
its moments (See Theorem \ref{em}). Thirdly, we establish a
super-process limit of such branching Markov processes where the
motion is Brownian (See Theorem \ref{super}).

The rest of the paper is organised as follows. In Section
\ref{sbps} we define the branching Markov process precisely and in
Section \ref{mainresult} we state the three main theorems of this
paper and make some remarks on various possible generalisations of
our results.

In Section \ref{rbp} we prove four propositions on age-dependent
Branching processes which are used in proving Theorem
\ref{oneparticle} (See Section \ref{pop}). In Section \ref{rbp} we
also show that the joint distribution of ancestoral times for a
sample of $k \geq 1$ individuals chosen at random from the
population at time $t$ converges as $ t \rightarrow \infty$ (See
Theorem \ref{cab}). This result is of independent interest and is
a key tool that is needed in proving Theorem \ref{em} (See
Section \ref{pemp}).

In Section \ref{psp}, we prove Theorem \ref{super}, the key idea
being to scale the age and motion parameters differently. Given
this, the proof uses standard techniques for such limits.  Theorem
\ref{oneparticle} is used in establishing the limiting log-Laplace
equation. Tightness of the underlying particle system is shown in
Proposition \ref{tight} and the result follows by the method
prescribed in \cite{D}.

\section{Statement of Results}

\subsection{The Model}
\label{sbps}

\mbox{\hspace{1in}}

Each particle in our system will have two parameters, age in
$\R_{+}$ and location in $\R$. We begin with the description of
the particle system.

\be
\item {\bf Lifetime Distribution $G(\cdot)$:} Let $G(\cdot)$ be a
cumulative distribution function on $[0,\infty)$, with $G(0) =
0.$ Let $\mu = {\int_0^\infty sdG(s)} < \infty.$

\item {\bf Offspring Distribution ${\bf p}$ :}
Let ${\bf p} \equiv \{p_k\}_{k \geq 0}$  be a probability distribution
such that $p_0 < 1$, $m = \sum_{k=0}^\infty k p_k =1$ and that
$\sigma^2 = \sum_{k=0}^\infty k^2 p_k - 1 < \infty$.

\item {\bf Motion Process $\eta(\cdot)$:}  Let $\eta(\cdot)$ be a $\R$
valued Markov process starting at $0$.
\ee

{\bf Branching Markov Process $(G, ${\bf p}$, \eta)$:} Suppose we
are given a realisation of an age-dependent branching process with
offspring distribution ${\bf p}$ and lifetime distribution $G$
(See Chapter IV of \cite{an} for a detailed description). We
construct a branching Markov process by allowing each individual
to execute an independent copy of $\eta$ during its lifetime
$\tau$ starting from where its parent died.

Let $N_t$ be the number of particles alive at time $t$ and  \eq{
\label{confi} \Ci_t = \{ (a^{i}_t, X^{i}_t) : i = 1,2, \ldots, N_t\}}
 denote the age and position configuration of all the individuals
alive at time $t.$ Since $m=1$ and $G(0)=0$, there is no
explosion in finite time (i.e. $P(N_t < \infty) =1$) and
consequently $\Ci_t$ is well defined  for each $0 \leq t < \infty$
(See \cite{an}).

Let $\Bi(\R_+)$ (and $\Bi(\R)$) be the Borel $\sigma$-algebra on
$\R_+$ (and $\R$).  Let $\mfs$ be the space of finite Borel
measures on $\R_{+} \times \R$ equipped with the weak topology.
Let $\mfas := \{ \nu \in \mfs : \nu = \sum_{i=1}^n \delta_{a_i,
x_i}(\cdot,\cdot), n \in \N, a_i \in \R_+, x_i \in \R \}.$ For
any set $A \in \Bi(\R_+)$ and $B \in \Bi(\R),$ let $Y_t(A \times
B) $ be the number of particles at time $t$ whose age is in $A$
and position is in $B$. As pointed out earlier, $m<\infty$, $G(0)
= 0$ implies that $Y_t \in \mfas$ for all $t > 0$ if $Y_0$ does
so. Fix a function $ \phi \in C^+_{b}(\R_{+} \times \R),$ (the
set of all bounded, continuous and positive functions from
$\R_+\times\R$ to $\R_+$), and define \ben \label{bd} \la Y_t,
\phi \ra = \int \phi \,dY_t = \sum_{i=1}^{N_t} \phi(a^{i}_t,
X^{i}_t).  \een

Since $\eta(\cdot)$ is a Markov process, it can be seen that
$\{Y_t: t \geq 0\}$ is a Markov process and we shall call $Y
\equiv \{Y_t : t \geq 0 \}$ {\it the $(G,{\bf p},\eta)$- branching
Markov process}.

Note that $\Ci_t$ determines $Y_t$ and conversely. The Laplace
functional of $Y_t,$ is given by
\eq{\label{lf} L_t \phi(a,x) := E_{a,x} [e^{- \la \phi, Y_t \ra}]
\equiv
 E[ e^{- \la \phi, Y_t \ra} \mid Y_0 = \delta_{a,x}].}

From the independence intrinsic in $\{Y_t: t \geq 0 \}$, we have:
\eq{ E_{\nu_1 + \nu_2}[e^{-\la \phi,Y_t \ra}] = (E_{\nu_1}[e^{-\la
\phi,Y_t \ra}]) (E_{\nu_2}[e^{-\la \phi,Y_t \ra }]),} for any
$\nu_i \in \mfas $ where $E_{\nu_i}[e^{-\la \phi,Y_t \ra }] :=
E[e^{-\la \phi,Y_t \ra }\mid Y_0 = \nu_i]$ for $i = 1,2$. This is
usually referred to as the branching property of $Y$ and can be
used to define the process $Y$ as the unique measure valued Markov
process with state space $\mfas$ satisfying $L_{t+s} \phi(a,x) =
L_t (L_s (\phi))(a,x)$ for all $t,s \geq 0.$

\subsection{The Results} \label{mainresult}

\mbox{\hspace{1in}}

In this section we describe the main results of the paper.  Let
$A_t$ be the event $\{ N_t > 0 \},$ where $N_t$ is the number of
particles alive at time $t$. As $p_0 < 1,$ $P(A_t) > 0$ for all $ 0
\leq t < \infty$ provided $P(N_0 = 0) \ne 1.$

\bt \label{oneparticle}  {\bf (Limiting behaviour of a randomly chosen
particle)} \\
On the event $A_t = \{N_t > 0 \}$, let $(a_t, X_t)$ be the age and
position of a randomly chosen particle from those alive at time
$t$. Assume that $\eta(\cdot)$ is such that for all $ 0 \leq t <
\infty$ 
\bea
\label{et1} &&  E(\eta(t)) = 0, v(t) \equiv E(\eta^2(t)) <
\infty, ~ \sup_{0\leq s \leq t}v(s) < \infty,  \\ && \mbox{ and } \psi
\equiv \int_0^\infty v(s)G(ds) < \infty. \nmr  
\eea
  Then, conditioned on
$A_t$, $(a_t,\f{X_t}{\sqrt{t}})$ converges as $ t \rightarrow
\infty$, to $(U,V)$ in distribution, where $U$ and $V$ are
Independent with $U$ a strictly positive absolutely continuous
random variable with density proportional to $(1 - G(\cdot))$ and
$V$ is normally distributed with mean $0$ and variance $
\f{\psi}{\mu}$. \et

Next consider the scaled empirical measure $\tilde{Y}_t \in
\mfas$ given by $\tilde{Y}_t(A \times B) = Y_t(A \times
\sqrt{t}B),$ $A \in \Bi(\R_+), B \in \Bi(\R).$

\bt \label{em} {\bf (Empirical Measure)}\\ Assume (\ref{et1}).
Then, conditioned on $A_t = \{ N_t > 0\},$ the random measures $\{
\f{\tilde{Y}_t}{N_t} \}$ converges as $t \rightarrow \infty$ in
distribution to a random measure $\nu,$ characterised by its
moment sequence $m_k(\phi) \equiv E[\nu(\phi)^k],$ for $\phi \in
C_b(\R_+ \times \R),$ $k \geq 1$. \et

An explicit formula for $m_k(\phi)$ is given in (\ref{mkpf})
below.

Our third result is on the super-process limit. We consider a
sequence of branching Markov processes $(G_n, {\bf p}_n,
\eta_n)_{\{n \geq 1\}}$ denoted by $\{Y_t^n: t \geq 0 \}_{\{n \geq
1\}}$ satisfying the following:

\be
\item[(a)] {\bf Initial measure:} For $n \geq 1$,
$Y^n_0 = \pi_{n \nu}$, where $\pi_{n \nu}$ is a Poisson random
measure with intensity $n \nu,$ for some $\nu = \alpha \times \mu
\in \mfs.$

\item[(b)] {\bf Lifetime $G^n(\cdot)$:}  For all $n \geq 1$,
$G^n$ is an exponential distribution with mean $\f{1}{\lambda}$

\item[(c)] {\bf Branching ${\bf p}_n,\cdot$:} For $n \geq 1$, Let
$F_n(u) = \sum_{k=0}^\infty p_{n,k} u^k$ be the generating function of
the offspring distribution ${\bf p}_n \equiv \{p_{n,k}\}_{k \geq
0}$. We shall assume that $F_n$ satisfies,
\eq{ \lim_{n \rightarrow \infty} \sup_{0 \leq u \leq N} \parallel
  n^2(F_n(1 - u/n) - (1 - u/n)) -  u^2 \parallel
\rightarrow 0, \label{asf}}
for all $N > 0.$

\item[(d)] {\bf Motion Process $\eta_n(\cdot)$:} For all $n \geq
1$,  \eq{\eta_n(t) = \f{1}{\sqrt{n}} \int_0^t \sigma(u )dB(u),
\qquad t \geq 0, \label{mtn} } where $\{B(t): t \geq 0\}$ is a
standard Brownian motion starting at $0$ and $\sigma: \R_+
\rightarrow \R$ is a continuous function such that $\int_0^\infty
\sigma^2(s) dG(s) < \infty $. It follows that for each $n \geq
1,$ $\eta_n$ satisfies (\ref{et1}).

\ee

\begin{defi}
Let $\Ei$ be an independent exponential random variable with mean
$\f{1}{\lambda}$, $ 0 < \lambda < \infty$. For $f \in C_l^+( \R_+
\times \R)$ let $U_t f (x) = E(f(\Ei, x + \sqrt{\lambda \psi}
B_t))$ where $\psi $ is defined in (\ref{et1}). For $t \geq 0$,
let $u_t(f)$ be the unique solution of the non linear integral
equation \eq{ u_t f(x) = U_t f(x) - \lambda \int_0^t U_{t-s}(u_s(
f)^2) (x)  ds. \label{R}}
Let $\{\Yi_t : t \geq 0\} $ be a $\mfs$ valued Markov process whose
Laplace functional is given by
\eq{ E_{\Ei \times \mu} [e^{-\la f, \Yi_t \ra}] = e^{-\la V_t f,
\mu \ra}, \label{log_Lap_super}}
where $f\in C_l^+(\R_{+} \times \R^d)$ (the set of all continuous
functions from $\R_+\times\R$ to $\R$ with finite limits as $
(a,x)\rightarrow \infty $) and $V_t(f)(x) \equiv u_t(f(x))$ for
$x \in \R$ (See \cite{D} for existence of $\Yi$ satisfying
(\ref{log_Lap_super}).
\end{defi}

Note that in the process $\{\Yi_t: t \geq 0 \}$ defined above, the
distribution of the age (i.e. the first coordinate) is
deterministic. The spatial evolution behaves like that of a
super-process where the motion of particles is like that of a
Brownian motion with variance equal to the average variance of
the age-dependent particle displacement over its lifetime. Also,
$u_s(f)$ in second term of (\ref{R}) is interpreted in the
natural way as a function on $\R_+ \times \R$ with $u_s(f)(a,x) =
u_s(f)(x)$ for all $a >0, x \in \R.$

\bt \label{super} {\bf (Age Structured Super-process)} \\ Let
$\epsilon > 0$. Let $\{Y_t^n: t \geq 0\}$ be the sequence of
branching Markov processes defined above(i.e.in (a), (b), (c),
(d)). Then as $n \rightarrow \infty$, $\{{\mathcal Y}_t^n \equiv
\f{1}{n} Y_{nt}^n, t \geq \epsilon \}$ converges weakly on the
Skorokhod space $\demas $ to $ \{\Yi_t : t \geq \epsilon\}$.  \et

\subsection{Remarks} \mbox{\hspace{1in}}

\medskip

(a) If $\eta(\cdot)$ is not Markov then $\tilde{C}_t = \{ a^i_t,
X^i_t , \tilde{\eta}_{t,i}\equiv\{ \eta_{t,i}(u) : 0 \leq u \leq
a^i_t\} : i = 1,2 \ldots, N_t \}$ is a Markov process where $\{
\tilde{\eta}_{t,i} (u) : 0 \leq u \leq a^i_t\} $ is the history of
$\eta(\cdot)$ of the individual $i$ during its lifetime.  Theorem
\ref{oneparticle} and Theorem \ref{em} extends to this case.

\medskip

(b) Most of the above results also carry over to the case when the
motion process is $\R^d$ valued ($ d \geq 1$) or is Polish space
valued and where the offspring distribution is age-dependent.

\medskip

(c) Theorem \ref{oneparticle} and Theorem\ref{em} can also be extended
to the case when $\eta(L_1)$, with $L_1 \eid G,$ is in the domain of
attraction of a stable law of index $0 < \al \leq 2.$

\medskip

(d) In Theorem \ref{super} the convergence should hold on $\dmas$ if
we take $\alpha$ in the sequence of branching Markov processes to be
$\Ei$ (i.e. Exponential with mean $\f{1}{\lambda}$).

(e) The super-process limit obtained in Theorem \ref{super} has been
considered in two special cases in the literature. One is in
\cite{bk} where an age-dependent Branching process is rescaled
(i.e. the particles do not perform any motion).  The other is in
\cite{DGL} where a general non-local super-process limit is
obtained when the offspring distribution is given by $p_1 = 1$.
In our results, to obtain a super-process limit the age-parameter
is scaled differently when compared to the motion parameter
giving us an age-structured super-process.

\medskip

(f) Limit theorems for critical branching Markov processes where the
motion depends on the age does not seem to have been considered in the
literature before.

\section{Results on Branching Processes} \label{rbp}

Let $\{ N_t : t \geq 0\}$ be an age-dependent branching process with
offspring distribution $\{p_k\}_{k \geq 0}$ and lifetime distribution
$G$ (see \cite{an} for detailed discussion).  Let $\{\zeta_k\}_{k\geq
0}$ be the embedded discrete time Galton-Watson branching process with
$\zeta_k$ being the size of the $k$th generation, $k \geq 0$. Let
$A_t$ be the event $\{N_t > 0\}$. On this event, choose an individual
uniformly from those alive at time $t$. Let $M_t$ be the generation number
and $a_t$ be the age of this individual.

\bp \label{lemma0} Let $A_t,a_t, M_t$ and $N_t$ be as above. Let
$\mu$ and $\sigma$ be as in Section~\ref{sbps}. Then \beas
 (a)&&\lim_{t \rightarrow \infty} t P(A_t) = \f{2 \mu}{\sigma^2}\\
 (b)&&  \mbox{For all } x >0 ,\, \, \lim_{t \rightarrow \infty}  P(\f{N_t}{t} > x|A_t) = e^{-\f{2 \mu x}{\sigma^2}},\\
 (c)&& \mbox{For all } \epsilon >0 ,\,\, \lim_{t \rightarrow \infty}  P( |\f{M_t}{t} - \f{1}{\mu} |  > \epsilon  |A_t) = 0\\
 (d) &&  \mbox{For all } x >0, \,\, \lim_{t \rightarrow \infty}  P(a_t \leq x | A_t) = \f{1}{\mu} \int_0^x (1-G(s))ds.
\eeas
\ep
{\em Proof :} For (a) and (b) see chapter 4 in \cite{an}.  For (c)  see  \cite{durrett} and for (d) see \cite{athreya}. \qed

\bp \label{llnh} (Law of large numbers) Let $\epsilon >0$ be
given. For the randomly chosen individual at time $t$, let $\{ L_{ti}
: 1 \leq i \leq M_t\}$, be the lifetimes of its ancestors.  Let $h:
[0,\infty) \rightarrow \R$ be Borel measurable and $E(\mid h(L_1)
\mid) < \infty$ with $L_1 \eid G.$ Then, as $t \rightarrow \infty$
$$P( | \f{1}{M_t} \sum_{i=1}^{M_t} h(L_{ti}) - E(h(L_1)) | > \epsilon | A_t) \rightarrow 0.$$
\ep

{\em Proof :} Let $\epsilon $ and $\epsilon_1 >0$ be given and let
$k_1(t) = t(\f{1}{\mu} -\epsilon)$ and $k_2(t) = t(\f{1}{\mu} + \epsilon).$ By
Proposition \ref{lemma0} there exists $\delta >0$, $\eta >0 $ and $t_0 > 0$
such that for all $t \geq t_0$, \eq{\label{1l}t P(N_t >0) > \delta
\mbox{ and } P(N_t \leq t \eta| A_t) < \epsilon_1 ;} \eq{\label{2l} P(
M_t \in [k_1(t), k_2(t)]^c |A_t) < \epsilon_1.} Also by the law of
large numbers for any $\{L_i\}_{i \geq 1}$ i.i.d. $G$ with $E |h(L_1)|
< \infty$ \eq{\label{3l} \lim_{k \rightarrow \infty} P (\sup_{j \geq
k}\f{1}{j} |\sum_{i=1}^{j} h(L_{i}) - E(h(L_1)) | > \epsilon) =
0.} Let $\{\zeta_{k}\}_{k \geq 0}$ be the embedded Galton-Watson process. For
each $t > 0$ and $k \geq 1$ let $\zeta_{kt}$ denote the number of
lines of descent in the $k$-th generation alive at time $t$ (i.e. the
successive life times $\{L_i\}_{i \geq 1}$ of the individuals in that
line of descent satisfying $\sum_{i=1}^k L_i \leq t \leq
\sum_{i=1}^{k+1}L_i$). Denote the lines of descent of these
individuals by $\{\zeta_{ktj}: 1 \leq j \leq \zeta_{kt}\}$. Call
$\zeta_{ktj}$ {\em bad} if \eq{
\label{bad} | \f{1}{k} \sum_{i=1}^{k} h(L_{ktji}) - E(h(L_1))) | >
\epsilon,} where $\{ L_{ktji} \}_{i \geq 1}$ are the successive
lifetimes in the line of descent $\zeta_{ktj}$ starting from the
ancestor. Let $\zeta_{kt,b}$ denote the cardinality of the set $\{
\zeta_{ktj} : 1 \leq j \leq \zeta_{kt} \mbox{ and } \zeta_{ktj} \mbox{
is bad} \}.$ Now, \bea \lefteqn{P( | \f{1}{M_t} \sum_{i=1}^{M_t}
h(L_{ti}) - E(h(L_1)) | > \epsilon | A_t)} \nmr\\&& \nmr\\ &=& P( \mbox{ The
chosen line of descent at time $t$ is {\em bad }} | A_t) \nonumber \\ && \nmr \\ &\leq&
P( \mbox{ The chosen line of descent at time $t$ is {\em bad}}, M_t \in
[k_1(t), k_2(t)] ) | A_t)\nmr \\ && + P( M_t \in [k_1(t)), k_2(t)]^c
|A_t) \nmr\\ &=& \f{1}{P(N_t > 0)} E(\f{
\sum_{j=k_1(t)}^{k_2(t)}\zeta_{jt,b}}{N_t}; A_t) + P( M_t \in
[k_1(t)), k_2(t)]^c |A_t)  \nmr\\ &=& \f{1}{P(N_t > 0)} E(\f{
  \sum_{j=k_1(t)}^{k_2(t)}\zeta_{jt,b}}{N_t}; N_t > t \eta)  + \nmr \\ &&    
+ \f{1}{P(N_t > 0)} E(\f{ \sum_{j=k_1(t)}^{k_2(t)}\zeta_{jt,b}}{N_t};
N_t \leq t \eta)  + P( M_t \in [k_1(t)),
k_2(t)]^c |A_t) \nmr\\  &\leq& \f{1}{P(N_t > 0)} E(
\f{\sum_{j=k_1(t)}^{k_2(t)}\zeta_{jt,b}}{t \eta}; N_t > t \eta) + \nmr \\ &&
\push  +\f{P( N_t \leq t \eta)}{P(N_t >0)} + P( M_t \in [k_1(t)), k_2(t)]^c
|A_t)\nmr \\ &=& \f{1}{t \eta P(N_t > 0)} \sum_{j=k_1(t)}^{k_2(t)}
E(\zeta_{jt,b} ) + \nmr \\ && + P( N_t \leq t \eta | N_t >0) + P( M_t \in [k_1(t)),
k_2(t)]^c |A_t) \nmr\\\label{is2} \eea For $t \geq t_0$ by (\ref{2l}) and
(\ref{3l}), the last two terms in (\ref{is2}) are less than
$\epsilon_1$. The first term is equal to \bea  \f{1}{t \eta
P(N_t > 0)} \sum_{j=k_1(t)}^{k_2(t)} E(\zeta_{jt,b} ) &=&\f{1}{t \eta P(N_t > 0)} \sum_{j=k_1(t)}^{k_2(t)}
E(\sum_{i=1}^{\zeta_j} 1_{\{\zeta_{jti} \mbox{ is bad.}
\}}) \nmr\eea \bea &=&\f{1}{t \eta P(N_t > 0)} \sum_{j=k_1(t)}^{k_2(t)}
E({\zeta_j}) \times \nmr \\&& \times P \left( \sum_{i=1}^j L_i \leq t < \sum_{i=1}^{j+1}L_i,
\f{1}{j} |\sum_{i=1}^{j} h(L_{i}) - E(h(L_1)) | > \epsilon \right),
\nmr \\
&&\nmr\\
&& \mbox{where the $\{L_i \}_{i \geq 1}$ are i.i.d. $G$.}\nmr
\eea
Using (\ref{1l}) and (since $m=1$) $E(\zeta_j) = E(\zeta_0)$ we can conclude that

\bea \lefteqn{ \f{1}{t \eta
P(N_t > 0)} \sum_{j=k_1(t)}^{k_2(t)} E(\zeta_{jt,b} ) } \nmr\\
&\leq& E(\zeta_0) \f{P (\sup_{j \geq k_1(t)}\f{1}{j}
|\sum_{i=1}^{j} h(L_{i}) - E(h(L_1)) | > \epsilon)}{t \eta P(N_t > 0)}
\nmr\\ &\leq& E(\zeta_0) \f{P (\sup_{j \geq k_1(t)}\f{1}{j}
|\sum_{i=1}^{j} h(L_{i}) - E(h(L_1)) | > \epsilon)}{\eta \delta}, \nmr \\
\label{is3} \eea which by (\ref{3l}) goes to zero.  So we have
shown that for $ t \geq t_0$,
$${P( | \f{1}{M_t} \sum_{i=1}^{M_t} h(L_{ti}) -
E(h(L_1)) | > \epsilon | A_t)} <  3 \epsilon_1. $$
Since $\epsilon_1 >0$ is arbitrary, the proof is complete.
\qed

\bp \label{clt} Assume (\ref{et1}) holds. Let $\{L_i\}_{i \geq
1}$ be i.i.d $G$ and $\{\eta_i\}_{i \geq 1}$ be i.i.d copies of
$\eta$ and independent of the $\{L_{i}\}_{i \geq 1}$. For $\theta
\in \R, t \geq 0$ define $\phi(\theta,t) = E e^{i \theta
\eta(t)}.$ Then there exists an event $D,$ with $P(D) =1$ and on
$D$ for all $\theta \in \R,$
$$ \prod_{j=1}^n \phi\left(\f{\theta}{\sqrt{n}}, L_j \right) \rightarrow e^{\f{-\theta^2\psi}{2}},
\qquad \mbox{ as } n \rightarrow \infty,$$

where $\psi$ is as in (\ref{et1}).\ep

{\em Proof:} Recall from (\ref{et1}) that $v(t) = E(\eta^2(t))$
for $ t \geq 0$. Consider $$X_{ni} =
\f{\eta_i(L_i)}{\sqrt{\sum_{j=1}^n v(L_j)}}\mbox{ for } 1 \leq i
\leq n$$ and ${\mathcal F} = \sigma(L_i : i \geq1).$ Given
${\mathcal F},$ $\{X_{ni} : 1 \leq i \leq n\}$ is a triangular
array of independent random variables such that for $1 \leq i
\leq n$, $E(X_{ni} | {\mathcal F}) = 0,$ $\sum_{i=1}^n E(X^2_{ni}
| {\mathcal F}) = 1.$

Let $\epsilon >0$ be given. Let $$L_n(\epsilon) = \sum_{i=1}^n E
\left( X_{ni}^2 \, ; X_{ni}^2 > \epsilon | {\mathcal F}
\right).$$ By the strong law of large numbers,
\ben\f{\sum_{j=1}^n v(L_j)}{n} \rightarrow \psi \push \mbox{ w.p.
1.} \label{llnv}\een

Let $D$ be the event on which (\ref{llnv}) holds.  Then on $D$ \beas
\limsup_{n \rightarrow \infty} L_n(\epsilon) &\leq& \limsup_{n
\rightarrow \infty} \f{\psi}{2n} \sum_{i=1}^n E( |\eta_i(L_i)|^2) : \mid
\eta_i(L_i)|^2 > \f{\epsilon n\psi}{2} | {\mathcal F})\\ &\leq&
\limsup_{k \rightarrow \infty} \f{\psi}{2} E(|\eta_1(L_1)|^2 : \mid
\eta_1(L_1) \mid^2 > k) \\ &=&0.  \eeas Thus the Linderberg-Feller Central
Limit Theorem  (see \cite{al}) implies, that on $D$, for all $\theta \in \R$
$$  \prod_{i=1}^n \phi\left(\f{\theta}{\sqrt{\sum_{j=1}^n
v(L_j)}}, L_j \right) = E(e^{i\theta \sum_{j=1}^n X_{nj}} | {\mathcal F}) \rightarrow  e^{\f{-\theta^2}{2}}.$$
Combining this with (\ref{llnv}) yields the result.
\qed

\bp  \label{p3.4} For the randomly chosen individual at time $t$,
let \\$\{ L_{ti}, \{\eta_{ti}(u) : 0 \leq u \leq L_{ti}\}:  1 \leq i
\leq M_t\}$, be the lifetimes and motion processes of its
ancestors. Let $Z_{t1} = {\f{1}{\sqrt{M_t}}\sum_{i=1}^{M_t}
\eta_{ti}(L_{ti})},$  and \\$\Li_t = \sigma \{M_t,  L_{ti} : 1 \leq i
\leq M_t\}.$ Then \ben E\left(  | E(e^{i \theta Z_{t1}} | \Li_t )
- e^{-\f{\theta^2 \psi}{2}}| \; | A_t \right)  \rightarrow 0  \een
\ep

{\em Proof:} Fix $\theta \in \R, \epsilon_1 >0$ and $\epsilon
>0$. Replace the definition of ``bad'' in (\ref{bad}) by
\eq{\label{badc} |\prod_{i=1}^{k} \phi(\f{\theta}{\sqrt{k}}, L_{ktji})
  - e^{-\f{\theta^2 \psi}{2}} | > \epsilon}

By Proposition \ref{clt} we have, \eq{\label{4l} \lim_{k \rightarrow
  \infty} P (\sup_{j \geq k}|\prod_{i=1}^{j}
  \phi(\f{\theta}{\sqrt{j}}, L_i) - e^{-\f{\theta^2 \psi}{2}} | >
  \epsilon) =0.}

 Using this in place of (\ref{3l}) and imitating the proof of
  Proposition \ref{llnh}, (since the details mirror that proof we
  avoid repeating them here), we obtain that for $t$ sufficiently
  large \eq{\label{clt1} P( |\prod_{i=1}^{M_t}
  \phi(\f{\theta}{\sqrt{M_t}}, L_{ti}) - e^{-\f{\theta^2 \psi}{2}} | >
  \epsilon_1 | A_t) < \epsilon .}  Now for all $ \theta \in \R$,
$$ E(e^{i \theta Z_{t1}} | \Li_t ) = \prod_{i=1}^{M_t} \phi(\f{\theta}{\sqrt{M_t}}, L_{ti}).$$ So,

\beas \lefteqn{\limsup_{t \rightarrow \infty} E(| E(e^{i\theta \f{1}{\sqrt{M_t}}
\sum_{i=1}^{M_t} \eta_i(L_{ti})}| \Li_t ) - e^{-\f{\theta^2 \psi}{2}}|  | A_t)}\\
&= & \limsup_{t \rightarrow \infty}E(|\prod_{i=1}^{M_t}
\phi(\f{\theta}{\sqrt{M_t}}, L_{ti}) - e^{-\f{\theta^2 \psi}{2}} |
| A_t)\\ &<& \epsilon_1 + 2  \limsup_{t \rightarrow \infty} P( |\prod_{i=1}^{M_t} \phi(\f{\theta}{\sqrt{M_t}}, L_{ti}) - e^{-\f{\theta^2 \psi}{2}} |
> \epsilon_1 | A_t)\\ &=& \epsilon_1 + 2 \epsilon.  \eeas
Since $\epsilon >0, \epsilon_1 >0$ are arbitrary we have the result.  \qed

The above four Propositions will be used in the proof of Theorem
\ref{oneparticle}. For the proof of Theorem \ref{em} we will need a
result on coalescing times of the lines of descent.

Fix $k \geq 2$. On the event $A_t = \{ N_t >0 \},$ pick $k$
individuals $C_1, C_2, \ldots, C_k$ from those alive at time $t$
by simple random sampling without replacement. For any two particles
$C_i, C_j$, let $\tau_{C_j,C_i,t}$ be the birth time of their most
recent common ancestor. Let $\tau_{k-1,t} = \sup \{
\tau_{C_j,C_i,t} : i \neq j, 1 \leq i,j \leq k \}$. Thus
$\tau_{k-1,t} $ is the first time there are $k-1$ ancestors of
the $k$ individuals $C_1, C_2, \ldots, C_k.$ More generally, for
$1\leq j \leq k-1$ let $\tau_{j,t}$ as the first time there are
$j$ ancestors of the $k$ individuals $C_1, C_2, \ldots C_k$.

\bt \label{cab}\mbox{}
\medskip

\be
\item For any $i,j$, $\lim_{t \rightarrow \infty } P(
\f{\tau_{C_i,C_j, t}}{t} \leq x | A_t) \equiv H(x)$ exists for all $x
\geq 0$ and $H(\cdot)$ is an absolutely continuous distribution
function on $[0,\infty] $
\item Conditioned on $A_t$ the vector
$\tilde{\tau}_t = \f{1}{t} (\tau_{j,t} : 1 \leq j \leq k-1)$ as $
t\rightarrow \infty$ converges in distribution to a random vector
$\tilde{T} = (T_1,\ldots, T_{k-1})$ with $0 < T_1<T_2< \ldots <
T_{k-1} <1$ and having an absolutely continuous distribution on
$[0,1]^{k-1}$. \ee \et

{\em Proof :} The proof of (i) and (ii) for cases $k=2,3$ is in
\cite{durrett}. The following is an outline of a proof of (ii) for the
case $k > 3$ (for a detailed proof see \cite{athreya}).

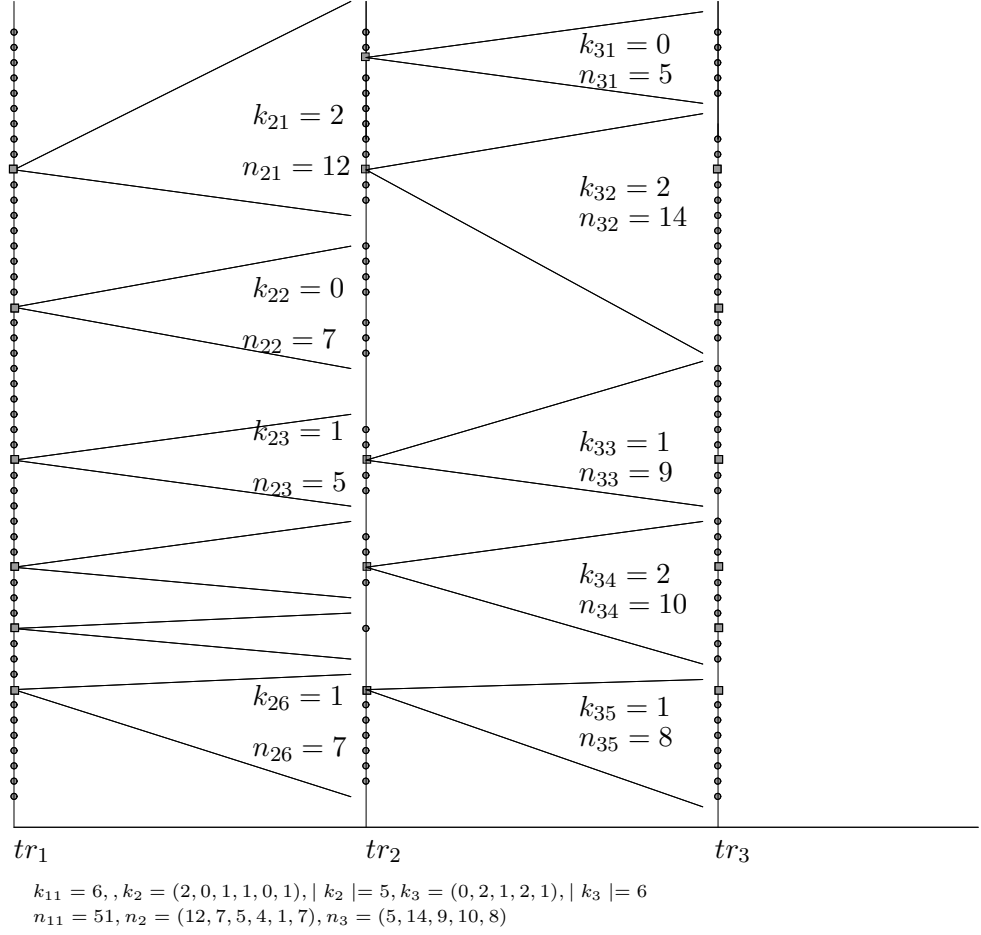
\begin{figure} \label{trp}
\setlength{\unitlength}{0.00053333in}
\begingroup\makeatletter\ifx\SetFigFont\undefined%
\gdef\SetFigFont#1#2#3#4#5{%
  \reset@font\fontsize{#1}{#2pt}%
  \fontfamily{#3}\fontseries{#4}\fontshape{#5}%
  \selectfont}%
\fi\endgroup%
{
\begin{picture}(9519,8764)(0,-10)
\texture{55888888 88555555 5522a222 a2555555 55888888 88555555 552a2a2a 2a555555
    55888888 88555555 55a222a2 22555555 55888888 88555555 552a2a2a 2a555555
    55888888 88555555 5522a222 a2555555 55888888 88555555 552a2a2a 2a555555
    55888888 88555555 55a222a2 22555555 55888888 88555555 552a2a2a 2a555555 }
\put(3507,8437){\shade\ellipse{60}{60}}
\put(3507,8437){\ellipse{60}{60}}
\put(3507,8287){\shade\ellipse{60}{60}}
\put(3507,8287){\ellipse{60}{60}}
\put(3507,7987){\shade\ellipse{60}{60}}
\put(3507,7987){\ellipse{60}{60}}
\put(3507,7837){\shade\ellipse{60}{60}}
\put(3507,7837){\ellipse{60}{60}}
\put(3507,7687){\shade\ellipse{60}{60}}
\put(3507,7687){\ellipse{60}{60}}
\put(3507,7537){\shade\ellipse{60}{60}}
\put(3507,7537){\ellipse{60}{60}}
\put(3507,7237){\shade\ellipse{60}{60}}
\put(3507,7237){\ellipse{60}{60}}
\put(3507,6937){\shade\ellipse{60}{60}}
\put(3507,6937){\ellipse{60}{60}}
\put(3507,6787){\shade\ellipse{60}{60}}
\put(3507,6787){\ellipse{60}{60}}

\put(3507,7387){\shade\ellipse{60}{60}}
\put(3507,7387){\ellipse{60}{60}}
\put(3507,6337){\shade\ellipse{60}{60}}
\put(3507,6337){\ellipse{60}{60}}
\put(3507,6187){\shade\ellipse{60}{60}}
\put(3507,6187){\ellipse{60}{60}}
\put(3507,6037){\shade\ellipse{60}{60}}
\put(3507,6037){\ellipse{60}{60}}
\put(3507,5887){\shade\ellipse{60}{60}}
\put(3507,5887){\ellipse{60}{60}}
\put(3507,5587){\shade\ellipse{60}{60}}
\put(3507,5587){\ellipse{60}{60}}
\put(3507,5437){\shade\ellipse{60}{60}}
\put(3507,5437){\ellipse{60}{60}}
\put(3507,5287){\shade\ellipse{60}{60}}
\put(3507,5287){\ellipse{60}{60}}

\put(3507,4387){\shade\ellipse{60}{60}}
\put(3507,4387){\ellipse{60}{60}}

\put(3507,4537){\shade\ellipse{60}{60}}
\put(3507,4537){\ellipse{60}{60}}
\put(3507,4087){\shade\ellipse{60}{60}}
\put(3507,4087){\ellipse{60}{60}}
\put(3507,3937){\shade\ellipse{60}{60}}
\put(3507,3937){\ellipse{60}{60}}
\put(3507,3487){\shade\ellipse{60}{60}}
\put(3507,3487){\ellipse{60}{60}}
\put(3507,3337){\shade\ellipse{60}{60}}
\put(3507,3337){\ellipse{60}{60}}
\put(3507,3037){\shade\ellipse{60}{60}}
\put(3507,3037){\ellipse{60}{60}}
\put(3507,1837){\shade\ellipse{60}{60}}
\put(3507,1837){\ellipse{60}{60}}
\put(3507,1687){\shade\ellipse{60}{60}}
\put(3507,1687){\ellipse{60}{60}}
\put(3507,1537){\shade\ellipse{60}{60}}
\put(3507,1537){\ellipse{60}{60}}
\put(3507,1387){\shade\ellipse{60}{60}}
\put(3507,1387){\ellipse{60}{60}}
\put(3507,1237){\shade\ellipse{60}{60}}
\put(3507,1237){\ellipse{60}{60}}
\put(3507,1087){\shade\ellipse{60}{60}}
\put(3507,1087){\ellipse{60}{60}}
\put(3507,2587){\shade\ellipse{60}{60}}
\put(3507,2587){\ellipse{60}{60}}
\put(6957,8437){\shade\ellipse{60}{60}}
\put(6957,8437){\ellipse{60}{60}}
\put(6957,8287){\shade\ellipse{60}{60}}
\put(6957,8287){\ellipse{60}{60}}
\put(6957,8137){\shade\ellipse{60}{60}}
\put(6957,8137){\ellipse{60}{60}}
\put(6957,7987){\shade\ellipse{60}{60}}
\put(6957,7987){\ellipse{60}{60}}
\put(6957,7837){\shade\ellipse{60}{60}}
\put(6957,7837){\ellipse{60}{60}}

\put(6957,7237){\shade\ellipse{60}{60}}
\put(6957,7237){\ellipse{60}{60}}
\put(6957,6937){\shade\ellipse{60}{60}}
\put(6957,6937){\ellipse{60}{60}}
\put(6957,6787){\shade\ellipse{60}{60}}
\put(6957,6787){\ellipse{60}{60}}
\put(6957,6487){\shade\ellipse{60}{60}}
\put(6957,6487){\ellipse{60}{60}}
\put(6957,6637){\shade\ellipse{60}{60}}
\put(6957,6637){\ellipse{60}{60}}
\put(6957,7387){\shade\ellipse{60}{60}}
\put(6957,7387){\ellipse{60}{60}}
\put(6957,6337){\shade\ellipse{60}{60}}
\put(6957,6337){\ellipse{60}{60}}
\put(6957,6187){\shade\ellipse{60}{60}}
\put(6957,6187){\ellipse{60}{60}}
\put(6957,6037){\shade\ellipse{60}{60}}
\put(6957,6037){\ellipse{60}{60}}
\put(6957,5887){\shade\ellipse{60}{60}}
\put(6957,5887){\ellipse{60}{60}}
\put(6957,5587){\shade\ellipse{60}{60}}
\put(6957,5587){\ellipse{60}{60}}
\put(6957,5437){\shade\ellipse{60}{60}}
\put(6957,5437){\ellipse{60}{60}}
\put(6957,5137){\shade\ellipse{60}{60}}
\put(6957,5137){\ellipse{60}{60}}

\put(6957,4987){\shade\ellipse{60}{60}}
\put(6957,4987){\ellipse{60}{60}}
\put(6957,4687){\shade\ellipse{60}{60}}
\put(6957,4687){\ellipse{60}{60}}
\put(6957,4387){\shade\ellipse{60}{60}}
\put(6957,4387){\ellipse{60}{60}}
\put(6957,4837){\shade\ellipse{60}{60}}
\put(6957,4837){\ellipse{60}{60}}
\put(6957,4537){\shade\ellipse{60}{60}}
\put(6957,4537){\ellipse{60}{60}}
\put(6957,4087){\shade\ellipse{60}{60}}
\put(6957,4087){\ellipse{60}{60}}
\put(6957,3937){\shade\ellipse{60}{60}}
\put(6957,3937){\ellipse{60}{60}}

\put(6957,3637){\shade\ellipse{60}{60}}
\put(6957,3637){\ellipse{60}{60}}
\put(6957,3487){\shade\ellipse{60}{60}}
\put(6957,3487){\ellipse{60}{60}}
\put(6957,3337){\shade\ellipse{60}{60}}
\put(6957,3337){\ellipse{60}{60}}
\put(6957,2887){\shade\ellipse{60}{60}}
\put(6957,2887){\ellipse{60}{60}}
\put(6957,2737){\shade\ellipse{60}{60}}
\put(6957,2737){\ellipse{60}{60}}
\put(6957,3037){\shade\ellipse{60}{60}}
\put(6957,3037){\ellipse{60}{60}}
\put(6957,2437){\shade\ellipse{60}{60}}
\put(6957,2437){\ellipse{60}{60}}
\put(6957,2287){\shade\ellipse{60}{60}}
\put(6957,2287){\ellipse{60}{60}}

\put(6957,1837){\shade\ellipse{60}{60}}
\put(6957,1837){\ellipse{60}{60}}
\put(6957,1687){\shade\ellipse{60}{60}}
\put(6957,1687){\ellipse{60}{60}}
\put(6957,1537){\shade\ellipse{60}{60}}
\put(6957,1537){\ellipse{60}{60}}
\put(6957,1387){\shade\ellipse{60}{60}}
\put(6957,1387){\ellipse{60}{60}}
\put(6957,1237){\shade\ellipse{60}{60}}
\put(6957,1237){\ellipse{60}{60}}
\put(6957,1087){\shade\ellipse{60}{60}}
\put(6957,1087){\ellipse{60}{60}}
\put(6957,937){\shade\ellipse{60}{60}}
\put(6957,937){\ellipse{60}{60}}
\put(57,8437){\shade\ellipse{60}{60}}
\put(57,8437){\ellipse{60}{60}}
\put(57,8287){\shade\ellipse{60}{60}}
\put(57,8287){\ellipse{60}{60}}
\put(57,8137){\shade\ellipse{60}{60}}
\put(57,8137){\ellipse{60}{60}}
\put(57,7987){\shade\ellipse{60}{60}}
\put(57,7987){\ellipse{60}{60}}
\put(57,7837){\shade\ellipse{60}{60}}
\put(57,7837){\ellipse{60}{60}}
\put(57,7687){\shade\ellipse{60}{60}}
\put(57,7687){\ellipse{60}{60}}
\put(57,7537){\shade\ellipse{60}{60}}
\put(57,7537){\ellipse{60}{60}}
\put(57,7237){\shade\ellipse{60}{60}}
\put(57,7237){\ellipse{60}{60}}
\put(57,6937){\shade\ellipse{60}{60}}
\put(57,6937){\ellipse{60}{60}}
\put(57,6787){\shade\ellipse{60}{60}}
\put(57,6787){\ellipse{60}{60}}
\put(57,6487){\shade\ellipse{60}{60}}
\put(57,6487){\ellipse{60}{60}}
\put(57,6637){\shade\ellipse{60}{60}}
\put(57,6637){\ellipse{60}{60}}
\put(57,7387){\shade\ellipse{60}{60}}
\put(57,7387){\ellipse{60}{60}}
\put(57,6337){\shade\ellipse{60}{60}}
\put(57,6337){\ellipse{60}{60}}
\put(57,6187){\shade\ellipse{60}{60}}
\put(57,6187){\ellipse{60}{60}}
\put(57,6037){\shade\ellipse{60}{60}}
\put(57,6037){\ellipse{60}{60}}
\put(57,5887){\shade\ellipse{60}{60}}
\put(57,5887){\ellipse{60}{60}}
\put(57,5587){\shade\ellipse{60}{60}}
\put(57,5587){\ellipse{60}{60}}
\put(57,5437){\shade\ellipse{60}{60}}
\put(57,5437){\ellipse{60}{60}}
\put(57,5137){\shade\ellipse{60}{60}}
\put(57,5137){\ellipse{60}{60}}
\put(57,5287){\shade\ellipse{60}{60}}
\put(57,5287){\ellipse{60}{60}}
\put(57,4987){\shade\ellipse{60}{60}}
\put(57,4987){\ellipse{60}{60}}
\put(57,4687){\shade\ellipse{60}{60}}
\put(57,4687){\ellipse{60}{60}}
\put(57,4387){\shade\ellipse{60}{60}}
\put(57,4387){\ellipse{60}{60}}
\put(57,4837){\shade\ellipse{60}{60}}
\put(57,4837){\ellipse{60}{60}}
\put(57,4537){\shade\ellipse{60}{60}}
\put(57,4537){\ellipse{60}{60}}
\put(57,4087){\shade\ellipse{60}{60}}
\put(57,4087){\ellipse{60}{60}}
\put(57,3937){\shade\ellipse{60}{60}}
\put(57,3937){\ellipse{60}{60}}
\put(57,3787){\shade\ellipse{60}{60}}
\put(57,3787){\ellipse{60}{60}}
\put(57,3637){\shade\ellipse{60}{60}}
\put(57,3637){\ellipse{60}{60}}
\put(57,3487){\shade\ellipse{60}{60}}
\put(57,3487){\ellipse{60}{60}}
\put(57,3337){\shade\ellipse{60}{60}}
\put(57,3337){\ellipse{60}{60}}
\put(57,2887){\shade\ellipse{60}{60}}
\put(57,2887){\ellipse{60}{60}}
\put(57,2737){\shade\ellipse{60}{60}}
\put(57,2737){\ellipse{60}{60}}
\put(57,3037){\shade\ellipse{60}{60}}
\put(57,3037){\ellipse{60}{60}}
\put(57,2437){\shade\ellipse{60}{60}}
\put(57,2437){\ellipse{60}{60}}
\put(57,2287){\shade\ellipse{60}{60}}
\put(57,2287){\ellipse{60}{60}}
\put(57,2137){\shade\ellipse{60}{60}}
\put(57,2137){\ellipse{60}{60}}
\put(57,1837){\shade\ellipse{60}{60}}
\put(57,1837){\ellipse{60}{60}}
\put(57,1687){\shade\ellipse{60}{60}}
\put(57,1687){\ellipse{60}{60}}
\put(57,1537){\shade\ellipse{60}{60}}
\put(57,1537){\ellipse{60}{60}}
\put(57,1387){\shade\ellipse{60}{60}}
\put(57,1387){\ellipse{60}{60}}
\put(57,1237){\shade\ellipse{60}{60}}
\put(57,1237){\ellipse{60}{60}}
\put(57,1087){\shade\ellipse{60}{60}}
\put(57,1087){\ellipse{60}{60}}
\put(57,937){\shade\ellipse{60}{60}}
\put(57,937){\ellipse{60}{60}}
\put(57,337){$tr_1$}
\put(3507,337){$tr_2$}
\put(6927,337){$tr_3$}
\path(57,7537)(57,8737)
\path(57,637)(9507,637)
\path(57,637)(9507,637)
\path(57,7087)(3357,8737)
\path(57,7087)(3357,8737)
\path(57,7087)(3357,6637)
\path(57,7087)(3357,6637)
\path(57,5737)(3357,6337)
\path(57,5737)(3357,6337)
\path(57,5737)(3357,5137)
\path(57,5737)(3357,5137)
\path(57,4237)(3357,4687)
\path(57,4237)(3357,4687)
\path(57,4237)(3357,3787)
\path(57,4237)(3357,3787)
\path(57,3187)(3357,3637)
\path(57,3187)(3357,3637)
\path(57,3187)(3357,2887)
\path(57,3187)(3357,2887)
\path(57,2587)(3357,2287)
\path(57,2587)(3357,2287)
\path(57,2587)(3357,2737)
\path(57,2587)(3357,2737)
\path(57,1987)(3357,2137)
\path(57,1987)(3357,2137)
\path(57,1987)(3357,937)
\path(57,1987)(3357,937)
\path(3507,8737)(3507,7387)
\path(3507,8737)(3507,7387)
\put(2300, 7532){ $k_{21} = 2$}
\put(2200, 7032){ $n_{21} = 12$}
\put(2300, 5832){ $k_{22} = 0$}
\put(2200, 5332){ $n_{22} = 7$}
\put(2300, 4432){ $k_{23} = 1$}
\put(2300, 3932){ $n_{23} = 5$}
\put(2300, 1832){ $k_{26} = 1$}
\put(2300, 1332){ $n_{26} = 7$}

\put(5500, 8232){ $k_{31} = 0$}
\put(5500, 7932){ $n_{31} = 5$}

\put(5500, 6832){ $k_{32} = 2$}
\put(5500, 6532){ $n_{32} = 14$}

\put(5500, 4332){ $k_{33} = 1$}
\put(5500, 4032){ $n_{33} = 9$}

\put(5500, 3032){ $k_{34} = 2$}
\put(5500, 2732){ $n_{34} = 10$}

\put(5500, 1732){ $k_{35} = 1$}
\put(5500, 1432){ $n_{35} = 8$}

\put(150, -200){{ \parbox{4in}{\tiny 
$k_{11} = 6,, k_2 = (2,0,1,1,0,1), \mid k_2 \mid= 5, k_3=(0,2,1,2,1), \mid k_3 \mid = 6$ \\
$n_{11} = 51, n_2 = (12,7,5,4,1,7),n_3=(5,14,9,10,8)$ }}}

\shade\path(3462,8232)(3537,8232)(3537,8157)
    (3462,8157)(3462,8232)
\shade\path(3462,7132)(3537,7132)(3537,7057)
    (3462,7057)(3462,7132)
\shade\path(3477,2017)(3552,2017)(3552,1942)
    (3477,1942)(3477,2017)
\path(3477,2017)(3552,2017)(3552,1942)
    (3477,1942)(3477,2017)
\shade\path(3477,4282)(3552,4282)(3552,4207)
    (3477,4207)(3477,4282)
\shade\path(3477,3232)(3552,3232)(3552,3157)
    (3477,3157)(3477,3232)
\path(3507,8187)(6807,8637)
\path(3507,8187)(6807,8637)
\path(3507,8187)(6807,7737)
\path(3507,8187)(6807,7737)
\path(3507,7537)(3507,637)
\path(3507,7087)(6807,7637)
\path(3507,7087)(6807,7637)
\path(3507,7087)(6807,5287)
\path(3507,7087)(6807,5287)
\path(3507,4237)(6807,5207)
\path(3507,4237)(6807,5207)
\path(3507,4237)(6807,3787)
\path(3507,4237)(6807,3787)
\path(3507,1987)(6807,2087)
\path(3507,1987)(6807,2087)
\path(3507,1987)(6807,837)
\path(3507,1987)(6807,837)
\path(6957,8737)(6957,7387)
\path(6957,8737)(6957,7387)
\path(6957,7537)(6957,637)
\shade\path(6912,7132)(6987,7132)(6987,7057)
    (6912,7057)(6912,7132)
\path(6912,7132)(6987,7132)(6987,7057)
    (6912,7057)(6912,7132)
\shade\path(6927,5767)(7002,5767)(7002,5692)
    (6927,5692)(6927,5767)
\path(6927,5767)(7002,5767)(7002,5692)
    (6927,5692)(6927,5767)
\shade\path(6927,4282)(7002,4282)(7002,4207)
    (6927,4207)(6927,4282)
\path(6927,4282)(7002,4282)(7002,4207)
    (6927,4207)(6927,4282)
\shade\path(6927,3232)(7002,3232)(7002,3157)
    (6927,3157)(6927,3232)
\path(6927,3232)(7002,3232)(7002,3157)
    (6927,3157)(6927,3232)
 \shade\path(6927,2017)(7002,2017)(7002,1942)
    (6927,1942)(6927,2017)
\path(6927,2017)(7002,2017)(7002,1942)
    (6927,1942)(6927,2017)
\shade\path(6927,2632)(7002,2632)(7002,2557)
    (6927,2557)(6927,2632)
\path(6927,2632)(7002,2632)(7002,2557)
    (6927,2557)(6927,2632)
\path(57,7537)(57,637)
\shade\path(12,7132)(87,7132)(87,7057)
    (12,7057)(12,7132)
\path(12,7132)(87,7132)(87,7057)
    (12,7057)(12,7132)
\shade\path(27,5767)(102,5767)(102,5692)
    (27,5692)(27,5767)
\path(27,5767)(102,5767)(102,5692)
    (27,5692)(27,5767)
\shade\path(27,4282)(102,4282)(102,4207)
    (27,4207)(27,4282)
\path(27,4282)(102,4282)(102,4207)
    (27,4207)(27,4282)
\shade\path(27,3232)(102,3232)(102,3157)
    (27,3157)(27,3232)
\path(27,3232)(102,3232)(102,3157)
    (27,3157)(27,3232)
\shade\path(27,2017)(102,2017)(102,1942)
    (27,1942)(27,2017)
\path(27,2017)(102,2017)(102,1942)
    (27,1942)(27,2017)
\shade\path(27,2632)(102,2632)(102,2557)
    (27,2557)(27,2632)
\path(27,2632)(102,2632)(102,2557)
    (27,2557)(27,2632)

\path(3507,3187)(6807,3637)
\path(3507,3187)(6807,3637)
\path(3507,3187)(6807,2237)
\path(3507,3187)(6807,2237)
\end{picture}
}

\caption{Tracking particles surviving at various times}
\end{figure}

Below, for $1 \leq j \leq k-1,$ $\tau_{j,t} $ will be denoted by
$\tau_j.$ It can be shown that it suffices to show that for any $
1 \leq i_1 < i_2 \ldots < i_p <k$ and $ 0 < r_1 < r_2< \ldots <
r_p <1$,
$$\lim_{t \rightarrow \infty } P( \f{\tau_{i_1}}{t} < r_1 <
\f{\tau_{i_2}}{t} < r_2 < \ldots < \f{\tau_{i_p}}{t} < r_p <
\f{\tau_{k-1}}{t} < r_{k-1} < 1 | A_t)$$ exists. We shall now
condition on the population size at time $t r_1$. Suppose that at
time $t r_1$ there are $n_{11}$ particles of which $k_{11}$ have
descendants that survive till time $tr_2$. For each $1 \leq j \leq
 k_{11},$ suppose there are  $n_{2j}$ descendants alive  at
time $tr_2$ and for each such $j$, let $k_{2j}$ out of the $n_{2j}$ have
descendants that survive till time $tr_3$. Let $k_2 = (k_{21},
\ldots, k_{2|k_1|})$ and $|k_2 | = \sum_{j=1}^{|k_1|} k_{2j}.$
Inductively at time $t r_i$, there are $n_{ij}$ descendants for
the $j$-th particle, $1 \leq j \leq | k_{i-1}|. $ For each such
$j$, let $k_{ij}$ out of $n_{ij}$ have descendants that survive
up till time $t r_{i+1}$ (See Figure \ref{trp} for an
illustration).

It will be useful to use the following notation: Let
$$n_{11} , k_{11} \in \N, k_{11} \leq n_{11},
| k_{1} \mid = k_{11},  n_{1} = (n_{11}).$$ For $i = 2, \ldots
i_p$ let $ (n_i , k_i) \in \N_i $, where $N_i \equiv \N^{\mid
k_{i-1}\mid } \times \N^{\mid k_{i-1} \mid}$
$$ k_{ij} \leq n_{ij}, \mid k_i \mid \equiv \sum_{j=1}^{|
k_{i-1}|} k_{ij}, {{n_i}\choose{k_i}} \equiv \prod_{j=1}^{|k_{i-1}|}
{{n_{ij}}\choose{k_{ij}}}.$$

Let $f_s = P(N_s > 0)$. Now,

\beas \lefteqn{P( \f{\tau_{i_1}}{t} < r_1 <
\f{\tau_{i_2}}{t} < r_2 < \ldots < \f{\tau_{i_p}}{t} < r_p <
\f{\tau_{k-1}}{t} < r_{k-1} < t | A_t)=} \\
&=& \f{f_{tr_1}}{f_t}\sum_{(n_i , k_i) \in
\N_i}\left ({{n_{11}}\choose{k_{11}}}(f_{tr_1})^{k_{11}}
(1-f_{tr_1})^{n_{11} -k_{11}} \right )
 \f{P(N_{tr_1} = n_1)}{f_{tr_1}} \times \\ &\times& \prod_{i=1}^{p+1} \prod_{j=1}^{|k_{i-1}|}
{{n_{ij}}\choose{k_{ij}}}(f_{tu_i})^{k_{ij}} (1-f_{u_i})^{n_{ij} -
k_{ij}} P(N^{j}_{tu_{i}} = n_{i,j}|
N^{j}_{t u_{i}} > 0)\times \\ && \times g({\bf k}) E \f{\prod_{j=1}^k X^j
}{S^k}, \eeas

with $u_i = r_{i+1} - r_{i}, i= 1,2,\ldots, p-1,$ $u_p = 1-r_{p}$,
$N^{j}_{tu_{i}}$ is number of particles alive at time $tu_i$ of the
age-dependent branching process starting with one particle namely $j$,
$g({\bf k}) = g(k_1, \ldots, k_p)$ is the proportion of configurations
that have the desired number of ancestors corresponding to the given
event, $X^j \eid N^{j}_{tu_p} | N^{j}_{t u_p} > 0$ and $S
= \sum_{j=1}^{| k_{p+1}|}X^j.$

\medskip

Let $q_i = \f{u_{i}}{u_{i+1}}$ for $1 \leq i \leq p-1$. Then
following \cite{durrett} and using Proposition \ref{lemma0} (i),
(ii) repeatedly we can show that $P( \f{\tau_{i_1}}{t} < r_1 <
\f{\tau_{i_2}}{t} < r_2 < \ldots < \f{\tau_{i_p}}{t} < r_p <
\f{\tau_{k-1}}{t} < r_{k-1} < t | A_t)$ converges to \beas
\lefteqn{\f{1}{q_1} \sum_{k_i \in \N^{\mid k_{i-1}\mid}} \int dx
e^{-x} (q_1 x)^{k_{11}} \f{1}{k_{11}!}e^{-x q_1}\times }\\ &&
\times \prod_{i=2}^{p+1}\prod_{j=1}^{|k_{i-1}|} \int dx e^{-x}
\f{(q_i x)^{k_ij}}{k_{ij}!} e^{-x q_i} )g({\bf k})\\ &&\times \int
\prod_{i=1}^{k+1}dx_i \left (\f{\prod_{i=1}^k
x_i}{(\sum_{i=1}^{k+1}x_i)^k}\right) e^{-\sum_{i=1}^{k+1}x_i}
\f{(x_{k+1})^{|k_{p+1}| -k}}{(|k_{p+1}| -k)!}\eeas \beas
 &=&\f{1}{q_1}\sum_{k_i
\in \N^{\mid k_{i-1}\mid}} \prod_{i=1}^{p+1}
\f{(q_i)^{|k_i|}}{(1+q_i)^{|k_i| - | k_{i-1}|}}g({\bf k}) \times \\ && \push \times\int
\prod_{i=1}^{k+1}dx_i \left (\f{\prod_{i=1}^k
x_i}{(\sum_{i=1}^{k+1}x_i)^k}\right) e^{-\sum_{i=1}^{k+1}x_i}
\f{(x_{k+1})^{|k_{p+1}| -k}}{(|k_{p+1}| -k)!}.  \eeas

Consequently, we have shown that the random vector $\tilde{\tau}_t$
converges in distribution to a random vector $\tilde{T}$.  From the
above limiting quantity, one can show that the $\tilde{T}$ has an
absolutely continuous distribution on $[0,1]^{k-1}$. See \cite{athreya}
for a detailed proof.

\qed

\section{Proof of Theorem \ref{oneparticle}} \label{pop}

For the individual chosen, let $(a_t, X_t)$ be the age and
position at time $t$.  As in Proposition \ref{p3.4},  let $\{
L_{ti}, \{\eta_{ti}(u), 0 \leq u \leq L_{ti}\} : 1 \leq i \leq
M_t\},$ be the lifetimes and the motion processes of the
ancestors of this individual  and $\{ \eta_{t(M_t+1)}(u) : 0 \leq u
\leq t-\sum_{i=1}^{M_t} L_{ti}\}$ be the motion this individual.
Let $\Li_t = \sigma(M_t, L_{ti}, 1 \leq i \leq M_t).$ It is
immediate from the construction of the process that:
$$ a_t = t - \sum_{i=1}^{M_t} L_{ti}, $$ whenever $M_t >0$ and is equal to $a+t$ otherwise; and  that $$ X_t =
X_0 + \sum_{i=1}^{M_t} \eta_{ti}(L_{ti}) + \eta_{t(M_t+1)}(a_t). $$

Rearranging the terms, we obtain
$$ (a_t, \f{X_t}{\sqrt{t}}) = (a_t, \sqrt{\f{1}{\mu}} Z_{t1}) + (0, \left (\sqrt{\f{M_t}{t}} -\sqrt{\f{1}{\mu}} \right)Z_{t2}) + (0,
 \f{X_0}{\sqrt{t}} + Z_{t2}) ,$$
where $Z_{t1} = \f{\sum_{i=1}^{M_t}
\eta_{ti}(L_{t_i})}{\sqrt{M_t}} and $ $Z_{t2} = \f{1}{\sqrt{t}}
\eta_{t(M_t +1)}(a_t)$.  Let $\epsilon >0$ be given.

\beas
P(|Z_{t2} | > \epsilon | A_t) &\leq & P(|Z_{t2} | > \epsilon, a_t  \leq  k | A_t) +   P(|Z_{t2} | > \epsilon, a_t > k | A_t)\\
&\leq &P(|Z_{t2} | > \epsilon, a_t  \leq  k | A_t) +   P(a_t > k | A_t)\\
&\leq& \f{E( |Z_{t2}|^2 I_{a_t  \leq  k} | A_t)}{\epsilon^2} + P(a_t > k | A_t)\\
\eeas

By Proposition \ref{lemma0} and the ensuing tightness, for any $\eta
>0$ there is a $k_\eta$
$$ P( a_t  >  k | A_t) < \f{\eta}{2}.$$
for all $ k \geq  k_\eta, t \geq 0$. Next,

\beas E( |Z_{t2}|^2 I_{a_t \leq k_\eta} | A_t) &=& E( I_{a_t\leq
k_\eta} E(|Z_{t2}|^2 | \Li_t) | A_t) \\ &=& E( I_{a_t\leq k_\eta}
\f{v(a_t)}{t} | A_t)\\ &\leq & \f{\sup_{u \leq k_\eta} v(u)}{t}.
\eeas Hence, \beas P(|Z_{t2} | > \epsilon | A_t) & \leq &
\f{\sup_{u \leq k_\eta} v(u)}{t \epsilon^2} + \f{\eta}{2} \eeas
Since $\epsilon >0 $ and $\eta>0$ are arbitrary this shows that
as $ t \rightarrow \infty$ \eq{\label{itl2} Z_{t2} | A_t \cid 0,}

Now, for $\lambda >0, \theta \in \R$, as $a_t$ is $\Li_t$
measurable we have \beas E(e^{-\lambda a_t} e^{-i
\f{\theta}{\sqrt{\mu}} Z_{t1}} | A_t) &= &   E( e^{-\lambda a_t}(
E(e^{-i \theta Z_{t1}} | \Li_t) - e^{-\f{\theta^2 \psi}{2 \mu}})
| A_t) +  \\ &&  + \, ^{-\f{\theta^2 \psi}{2 \mu}} E(e^{-\lambda a_t} | A_t)
 \eeas

Proposition \ref{clt} shows that the first term above converges to
zero and using Proposition \ref{lemma0} we can conclude that as $
t \rightarrow \infty$  \ben \label{itl3}
(a_t,\f{1}{\sqrt{\mu}}Z_{t1}) | A_t \cid (U,V) \een As $X_0$ is a
constant, by Proposition \ref{lemma0} (c), (\ref{itl3}),
(\ref{itl2}) and Slutsky's Theorem, the proof is complete.  \qed

\section{Proof of Theorem \ref{em}} \label{pemp}

Let $\phi \in C_b(\R \times \R_+)$. We shall show, for each $k \geq
1$, that the moment-functions of $E(\f{<\tilde{Y}_t, \phi>^k}{N^k_t} |
A_t)$ converges as $t \rightarrow \infty$. Then by Theorem 16.16 in
\cite{ka} the result follows.

 The case $k=1$ follows from Theorem \ref{oneparticle} and the bounded
 convergence theorem. We shall next consider the case $k=2.$ Pick two
 individuals $C_1,C_2$ at random (i.e. by simple random sampling
 without replacement) from those alive at time $t$. Let the age and
 position of the two individuals be denoted by $(a^i_t,X^i_t), i
 =1,2.$ Let $\tau_t = \tau_{C_1,C_2,t}$ be the birth time of their
 common ancestor, say $D$, whose position we denote by
 $\tilde{X}_{\tau_t}$. Let the net displacement of $C_1$ and $C_2$
 from $D$ be denoted by $X^i_{t-\tau_t}, i=1,2$ respectively. Then
 $X^i_t = \tilde{X}_{\tau_t} + X^{i}_{t -\tau_t}, i =1,2$.

 Next, conditioned on this history up to the birth of $D (\equiv
\Gi_t)$, the random variables $(a^i_t, X^i_{t-\tau_t}), i =1,2$
are independent. By Proposition \ref{cab} (i) conditioned on
$A_t$, $\f{\tau_t}{t}$ converges in distribution to an absolutely
continuous random variable $T$ (say) in $[0,1]$. Also by Theorem
\ref{oneparticle} conditioned on $\Gi_t$ and $A_t$, $\{
(a^i_t,\f{X^i_{t-\tau_t}}{\sqrt{t-\tau_t}}), i = 1,2\} $
converges in distribution to $\{(U_i,V_i), i =1,2\}$ which are
i.i.d. with distribution $(U,V)$ as in Theorem \ref{oneparticle}.
Also $\f{\tilde{X}_{\tau_t}}{\sqrt{\tau_t}}$ conditioned on
$A_{\tau_t}$ converges in distribution to a random variable $S$
distributed as $V$.

Combining these one can conclude that $\{ (a^i_t,
\f{X^i_t}{\sqrt{t}}), i =1,2\}$ conditioned on $A_t$ converges in
distribution to $\{(U_i, \sqrt{T}S + \sqrt{(1-T)}V_i), i =1,2\}$ where $U_1, U_2, S,
V_1, V_2$ are all independent. Thus for any $\phi \in C_b (\R_+ \times
\R)$ we have, by the bounded convergence theorem,

\eq{ \lim_{t \rightarrow \infty} E( \prod_{i=1}^2 \phi (a^i_t, \f{X^i_t}{\sqrt{t}}) | A_t)  = E \prod_{i=1}^2\phi(U_i, \sqrt{T}S+ \sqrt{(1-T)}V_i) \equiv m_2(\phi) \mbox{ (say)}}

Now,
\beas E ( \left(\f{\tilde{Y}_t(\phi )}{N_t} \right)^2 | A_t ) &=& E(\f{(\phi(a_t, \f{X_t}{\sqrt{t}}) )^2 }{N_t}| A_t)\\ && +  E( \prod_{i=1}^2 \phi (a^i_t, \f{X^i_t}{\sqrt{t}})  \f{N_t (N_t-1)}{N_t^2} | A_t) \eeas

Using Proposition \ref{lemma0} (b) and the fact that $\phi$ is bounded
we have \\$\lim_{t \rightarrow \infty} E ( (\f{\tilde{Y}_t(\phi)}{N_t}
)^2 | A_t ) $ exists in $(0,\infty)$ and equals $m_2(\phi)$. The case
$k > 2$ can be proved in a similar manner but we use Theorem \ref{cab} (ii)
as outlined below. First we observe that as $\phi$ is bounded,
$$ E \left (\f{<\tilde{Y}_t, \phi>^k}{N^k_t} | A_t \right )   + 
=  \sum_{{\bf i}} E h(N_t,k)   \left(\prod_{j=1}^k \phi(a^{i_j}_t, \f{X^{i_j}_t}{\sqrt{t}}) | A_t \right)+ g( \phi,\Ci_t, N_t),$$

where $ h(N_t,k) \rightarrow 1$ and $g( \phi,\Ci_t, N_t) \rightarrow 0$ as $ t \rightarrow \infty$; and  ${\bf i} = \{i_1,i_2,\ldots, i_k\}$ is the index of $k$particles sampled without  replacement from $\Ci_t$ (see
(\ref{confi})). Consider one such sample, and re-trace the genealogical
tree $\Ti_{\bf i} \in \Ti(k)$,($\Ti(k)$ is the collection of all
possible trees with $k$ leaves given by ${\bf i}$), until their most
common ancestor. For any leaf $i_j$ in $\Ti_{\bf i}$, let $1 =
n(i_j,1) < n(i_j,2) < \cdots < n(i_j,N_{i_j})$ be the labels of the
internal nodes on the path from leaf $i_j$ to the root. We list the
ancestoral times on this by $\{\tau_{1}, \tau_{n(i_j,1)}, \ldots,
\tau_{n(i_j,N_{i_j})}.$ Finally we denote the net displacement of the
ancestors in the time intervals
$$[0,\tau_1], [\tau_{1},\tau_{n(i_j,2)} ], \ldots, [\tau_{n(i_j,N_{i_j}-1)},
  \tau_{n(i_j,N_{i_j} )}], [\tau_{n(i_j,N_{i_j})}, t]$$ by
 \[\tilde{\eta}^1_{i_j}({\tau_1}),
\tilde{\eta}_{i_j}^2(\tau_{n(i_j,2)},\tau_1), \ldots,
\tilde{\eta}_{i_j}^{N_{i_j}}(\tau_{n(i_j,N_{i_j}
)},\tau_{n(i_j,N_{i_j}-1)
}),\tilde{\eta}_{i_j}^\prime(t,\tau_{n(i_j,N_{i_j})}).\]
Given the above notation we have:
\beas
{ E \left (\prod_{j=1}^k \phi(a^{i_j}_t, \f{X^{i_j}_t}{\sqrt{t}}) | A_t \right) =} E \left( \sum_{ T \in \Ti_{{\bf i }}}  \prod_{j=1}^k  f(\phi,j,t) | A_t \right), \eeas
where $$f(\phi,j,t) = \phi(a^{i_j}_t, \f{1}{\sqrt{t}}
(\tilde{\eta}_{i_j}^1({\tau_1}) + \sum_{m=2}^{N_{i_j}}
\tilde{\eta}_{i_j}^m(\tau_{n(i_j,m)},\tau_{n(i_j,m-1)}) +
\tilde{\eta}_{i_j}^\prime(t,\tau_{n(i_j,N_{i_j})}).$$
Now by Theorem \ref{cab}, $$\f{(\tau_1,\tau_{n(i_j,2)}, \ldots,
\tau_{n(i_j,N_{i_j})})}{\sqrt{t}} |A_t \cid (T_1,T_{n(i_j,2)}, \ldots,
T_{n(i_j,N_{i_j})}).$$ So by  Theorem \ref{oneparticle}
{ 
\bea \lim_{t \rightarrow \infty}  E ( \left(\f{\tilde{Y}_t(\phi )}{N_t} \right)^2 | A_t ) = E \left (\sum_{{\bf i}} \sum_{ T \in \Ti_{{\bf i }}} \prod_{j=1}^k g(\phi,j,t) | A_t  \right)  \nonumber \equiv m_k(\phi) \nmr
\\\label{mkpf} \eea }where \beas
\lefteqn{g(\phi,j,t) = }\\ &=&\phi\left (U, S \sqrt{T_1} + \sum_{m=2}^{N_{i_j}} Z_{i_j}^m
\sqrt{T_{n(i_j,m)}- T_{n(i_j,m-1)}} + Z_{i_j}^\prime \sqrt{1-
T_{n(i_j,N_{i_j})}} \right )\eeas with $S, Z^\prime_{i_j}, Z^m_{i_j},$ $m=2,\ldots
,N_{i_j}, $ are i.i.d.$V$, $U$ is an independent random variable given
in Theorem \ref{oneparticle} and $T_i$'s are as in Theorem \ref{cab}
(ii). Since $\phi$ is bounded, the sequence $\{ m_k(\phi) \equiv
\lim_{t \rightarrow \infty } E(\f{<\tilde{Y}_t, \phi>^k}{N^k_t} )\}$
is necessarily a moment sequence of a probability distribution on
$\R$. This being true for each $\phi$, by Theorem 16.16 in \cite{ka}
we are done.  \qed

\section{ Proof of Theorem \ref{super}} \label{psp}

Let $Z$ be the Branching Markov process $Y$ described earlier, with
lifetime $G$ exponential with mean $\lambda$, $p_1 =1$ and $\eta \eid
\eta_1$ (see (\ref{mtn})).  Then it is easy to see that for any
bounded continuous function, $S_t\phi(a,x) = E_{(a,x)}< Z_t, \phi> =
E_{(a,x)} \phi(a_t,X_t)$ satisfies the following equation: \ben
S_t\phi(a,x) = e^{-\lambda t} W_t\phi(a,x) + \int_0^t ds \lambda
e^{-\lambda s} W_{s}(S_{t-s}(\phi)(0,\cdot))(a,x), \een

where $W_t$ is the semi-group associated to $\eta_1.$ Let $\Li$ be the
generator of $\eta_1$. Making a change of variable $s
\rightarrow t-s$ in the second term of the above and then
differentiating it with respect to $t$, we have \bea
\lefteqn{\f{d}{dt}S_t(\phi)(a,x) = -\lambda e^{-\lambda t} W_t\phi(a,x) +
e^{-\lambda t}\Li W_t \phi(a,x) + \lambda S_{t}(\phi)(0,x)} \nonumber
\\&& \push +\int_0^t ds \lambda(-\lambda e^{-\lambda (t-s)})
W_{t-s}(S_{s}(\phi)(0,\cdot))(a,x) \nonumber \\&& \push +\int_0^t
ds\lambda e^{-\lambda (t-s)} \Li
W_{t-s}(S_{s}(\phi)(0,\cdot))(a,x)\nonumber \\ &=& \lambda
S_{t}(\phi)(0,x) \nmr \\ &&+ (\Li -\lambda )\left [e^{-\lambda t} W_t\phi(a,x) +
\int_0^t ds \lambda e^{-\lambda (t-s)}
W_{t-s}(S_{s}(\phi)(0,\cdot))(a,x) \right ]\nonumber \\ &=& \lambda
S_{t}(\phi)(0,x) + (\Li-\lambda) S_t(\phi)(a,x) \nonumber \\ &=&
\f{\partial S_t\phi}{\partial a}(a,x) + \f{\sigma^2(a)}{2}{\Delta
S_t\phi}(a,x) + \lambda (S_{t}(\phi)(0,x) - S_t(\phi)(a,x)) \nmr, \eea

For each $n \geq 1$ define (another semigroup) $R^n_t\phi(a,x) =
E_{a,0}(\phi(a_t, x + \f{X_t}{\sqrt{n}}).$ Now note that, \beas
R^n_t\phi(a,x) &=& E_{a,0}(\phi(a_t, x + \f{X_t}{\sqrt{n}})\\
               &=& E_{a,\sqrt{n}x}(\phi(a_t, \f{X_t}{\sqrt{n}}) \\
               &=& S_{t}\phi_n(a,\sqrt{n}x),
\eeas where $\phi_n(a,x) = \phi(a, \f{x}{\sqrt{n}}).$
Differentiating w.r.t. $t$, we have that the generator of $R^n_t$
is

\eq{ \Ri^n\phi(a,x) = \f{\partial \phi}{\partial a}(a,x) +
\f{\sigma^2(a)}{2n}{\Delta\phi}(a,x) + \lambda (\phi(0,x) -
\phi(a,x)).}

\bp \label{bcp} Let $\epsilon >0$ and $t \geq \epsilon$. Let $\phi
\in C_l^+ (\R_+ \times \R^d)$. Then, \eq{\sup_{(a,x) \in \R_+
\times \R} \mid R^n_{nt}(\phi)(a,x)-
  U_t(\phi)(x)\mid \rightarrow 0.} \ep

{\em Proof:} Let $t\geq \epsilon$. Applying Theorem \ref{oneparticle} to the
process $Z$, we have $(a_{nt}, \f{X_{nt}}{\sqrt{n}}) \cid (U,V)$. The
proposition is then immediate from the bounded convergence theorem and
the fact that $\phi \in C_l^+ (\R_+ \times\R)$
 \qed

\bp
 \label{logl} Let $\pi_{n\nu}$ be a Poisson random measure with intensity $n \nu$ and $ t \geq 0$. The  log-Laplace functional of ${\mathcal Y}_t^n$,
\eq{\label{scll} E_{\pi_{n \nu}} [e^{-\la \phi, {\mathcal Y}_t^n
\ra}] = e^{- \la u_t^n \phi, \nu \ra}, }
where \eq{ \label{neq} u^n_t \phi(a,x) =
R^n_{nt}n(1-e^{-\f{\phi}{n}})(a,x) - \lambda \int_0^t ds R^n_{n(t
-s)}(n^2 \Psi_n(\f{u^n_s \phi}{n}))(a,x),} where
\[ \Psi_n (\phi)(a,x)  :=  \left[  F_n(1-\phi(0,x)) -(1-\phi(0,x)) \right]. \]
\ep

{\em Proof:} For any $n\in \N$, let $Y^n_t$ be the sequence of
branching Markov processes defined in Section \ref{mainresult}. It can 
be shown that its log-Laplace functional $L^{n}_t$ satisfies,
\ben L^n_{nt} \phi (a,x) =  e^{-\lambda n
t}W^n_{nt}[e^{-\phi}](a,x) +  \int_0^{nt} ds \lambda e^{-\lambda
s}W^n_{s} \left [ F_n(L^{n}_{nt-s}\phi(0,\cdot)) \right](a,x) ds,
\een
where $ t \geq 0$ and $W^n_t$ is the semigroup associated with $\eta_n$.
Using the fact that $e^{-\lambda u} = 1 - \int_0^u ds \lambda
e^{-\lambda s}$ for all $u \geq 0$ and a routine simplification, as done in
\cite{dy}, will imply that
\eq{ L^n_{nt} \phi(a,x) = W^n_{nt}[e^{-\phi}](a,x) + \lambda
\int_0^{nt} W^n_{nt -s}( F_n(L^n_{s}\phi(0,\cdot)) - L^n_{s}\phi
)(a,x) ds}
Therefore $v^n_{nt}(\phi)(a,x) = 1- L^n_t \phi(a,x),$ satisfies,
\eq{\label{ev} v^n_{nt} \phi (a,x) = W^n_{nt} (1 -
e^{-\phi})(a,x) +  \int_0^{nt} ds  W^n_{nt-s} ( (1 - v^n_{s}\phi)
- F_n(1-v^n_{s})\phi(0,\cdot))   )(a,x) \lambda  ds.}
Let $\Li^n$ be the generator of $\eta_n$. Then for $0 \leq s < t$\beas
\lefteqn{\f{d}{ds} R^n_{n(t-s)}(v^n_{ns}(\phi))(a,x)= }\\ &=& -(n\Ri^n)R^n_{n(t-s)}\left( v^n_{ns} (\phi) \right)(a,x) + R^n_{n(t-s)} \left( \frac{\partial}{\partial s} v^n_{ns} (\phi) \right)(a,x)\\
&=& -(n\Ri^n)R^n_{n(t-s)}\left( v^n_{ns} (\phi) \right)(a,x) \\ && + R^n_{n(t-s)} \left( n \Li^n W^n_{ns}(1-e^{-\phi}) + n\lambda ((1 - v^n_{ns}\phi) - F_n(1-v^n_{ns})\phi(0,\cdot))   )(a,x) \right )\\ && +  R^n_{n(t-s)}\left (\int_0^{ns} dr n\Li^n(W^n_{ns-r}((1-v^n_r(\phi)) -F_n(1-v^n_r)\phi(0,\cdot)) \right)(a,x)\\
&=& R^n_{n(t-s)}n\left(  -\lambda (v^n_{ns}(\phi)(0,\cdot) -v^n_{ns}(\phi)) +  \lambda ((1 - v^n_{ns}\phi) - F_n(1-v^n_{ns})\phi(0,\cdot))   )(a,x) \right )\\
&=& -R^n_{n(t-s)}(n\Psi_n(v^n_{ns}\phi))(a,x), \eeas

Integrating both sides with respect to $s$ from $0$ to $t$, we obtain
that \ben\label{kv} v^n_{nt}(\phi)(a,x) = R^n_{nt}(1-e^{-\phi})(a,x)
-\int_0^{t} ds R^n_{n(t-s)}(n\Psi_n(v^n_{ns}\phi))(a,x). \een If
$\pi_{n \nu}$ is a Poisson random measure with intensity $n\nu$, then
\beas E_{\pi_{n \nu}} [e^{-\la \phi, {\mathcal Y}^n_t \ra}] =
E_{\pi_{n \nu}} [e^{-\la \frac{\phi}{n}, {Y}^n_{nt} \ra}]= e^{\la
L^n_t(\f{\phi}{n})-1,n\nu \ra}= e^{ -\la nv^n_t(\f{\phi}{n}),\nu
\ra}. \eeas Therefore if we set $u^n_t (\phi) \equiv nv^n_{nt}(\f{\phi}{n})$.
From (\ref{kv}), it is easy to see that $u^n_t(\phi)$ satisfies (\ref{scll}). \qed

\medskip

For any $f: \R_+ \times \R \rightarrow
\R,$ , we let $\parallel f \parallel_\infty = \sup_{(a,x) \in \R_+
\times \R} \mid f(a,x)\mid.$ With a little abuse of notation we shall
let $\parallel f \parallel_\infty = \sup_{x \in \R} \mid f(x)\mid$
when $f: \R \rightarrow \R$ as well.
\medskip

\bp \label{ucc} Let $ \epsilon >0$.  $\phi \in
C_l^+(\R_+\times\R^d)$ and $u^n_t(\phi)$ be as in Proposition
\ref{logl} and $u_t(\phi)$ be as in Theorem \ref{super}. Then for
$t \geq \epsilon$, \ben
 \sup_{(a,x)
\in \R_+ \times \R} \mid u^n_t(\phi)(a,x) - u_t(\phi)(x)\mid
\rightarrow 0 \een \ep

{\em Proof:} For any real $u \in \R$, define, $\eps_n(u) =
\lambda n^2(F_n(1-\f{u}{n}) - (1-\f{u}{n})) - u^2. $ So, \beas
\lefteqn{u^n_t(\phi)(a,x) = R^n_{nt}n(1-e^{-\f{\phi}{n}})(a,x) -\lambda
\int_0^t ds R^n_{n(t -s)}(n^2 \Psi_n(\f{u^n_s \phi}{n}))(a,x)}\\
&=&R^n_{nt}n(1-e^{-\f{\phi}{n}})(a,x) - \int_0^t ds R^n_{n(t
-s)}(\eps_n(u^n_s(\phi(0\cdot))))(a,x)\\ && \push - \lambda\int_0^t ds
R^n_{n(t -s)} (u^n_{s}\phi(0,\cdot)^2)(a,x) \\ \eeas Now
 \beas
\lefteqn{u^n_t(\phi)(a,x) - u_t(\phi)(x) =}\\&=&
R^n_{nt}(n(1-e^{-\f{\phi}{n}}))(a,x) -U_t(\phi)(x)\\ & - &
\int_0^t ds R^n_{n(t -s)}(\eps_n(u^n_s(\phi(0\cdot))))(a,x) \\
&+& \lambda \int_0^t ds \left(  U_{(t -s)}( (u_{s}\phi)^2)(a,x)- R^n_{n(t -s)}
(u^n_{s}\phi(0,\cdot)^2)(a,x) \right) \eeas \beas &=&
R^n_{nt}(n(1-e^{-\f{\phi}{n}}))(a,x) -U_t(\phi)(x) - \int_0^t ds
R^n_{n(t -s)}(\eps_n(u^n_s(\phi(0,\cdot))))(a,x) \\&& + \lambda \int_0^t
ds R^n_{n(t -s)}(( u_{s}\phi)^2-u^n_{s}\phi(0,\cdot)^2)(a,x) \\
&& \push + \lambda \int_0^t ds \left (U_{t-s}(u_{s}\phi)^2)(x) -R^n_{n(t -s)}
(u_{s}\phi)^2)(a,x) \right ) \\ \eeas Observe that, $R^n_{\cdot}$ is a
contraction, $\parallel u^n_\cdot(\phi)\parallel_\infty \leq
\parallel \phi \parallel_\infty$ and $\parallel u_\cdot(\phi)
\parallel_\infty \leq \parallel \phi \parallel_\infty$ for $\phi
\in C_l(\R_+\times \R).$ Therefore, we have \beas \parallel
u^n_t(\phi) - u_t(\phi)
\parallel_{\infty} &\leq &\parallel R^n_{nt}(n(1-e^{-\f{\phi}{n}}))
-U_t(\phi) \parallel_{\infty} + t \parallel
\epsilon_n(u^n_s(\phi(0,\cdot))\parallel_\infty \\ &&+ 2\lambda  \parallel
\phi \parallel_\infty\int_0^t ds \parallel u^n_s (\phi) -
u_s(\phi)
\parallel_\infty \\ &&+ \lambda \int_0^t ds \parallel (U_{t-s} -R^n_{n(t
-s)})(u_{s}\phi)^2 \parallel_\infty. \eeas  For $\phi \in
C_l(\R_+\times\R^d),$ note that, $U_t$, is a strongly continuous
semi-group implies that $u_s(\phi)$ is a uniformly continuous
function. So using Proposition \ref{ucc} the first term and the
last term go to zero.  By our assumption on $F$, $\parallel
\epsilon_n(u^n_s(\phi(0,\cdot))\parallel_\infty $ will go to zero
as $n \rightarrow \infty.$ Now using the standard Gronwall
argument we have the result. \qed

\bp \label{tight} Let $\epsilon >0$. The processes ${\mathcal Y}^n_\cdot$ are tight in $\demas$. \ep

{\bf Proof} By Theorem 3.7.1 and Theorem 3.6.5 (Aldous
Criterion) in \cite{D} , it is enough to show \eq{\label{tight0}\la {\mathcal Y}^n_{\tau_n +
\delta_n}, \phi\ra - \la {\mathcal Y}^n_{\tau_n}, \phi\ra \cid 0,}where $\phi \in C^+_l(\R_+\times \R)$,  $\delta_n$ is a sequence of positive numbers that converge
to $0$ and $\tau_n$ is any stop time of the process ${\mathcal Y}^n$
with respect to the canonical filtration, satisfying $0 <
\epsilon \leq \tau_n \leq T$ for some $T < \infty$.

First we note that, as $\la {\mathcal Y}^n_t, 1 \ra$ is a martingale,
for $\gamma >0$ by Chebyschev's inequality and Doob's maximal inequality we have\eq{  \label{tight1} P( \la {\mathcal
Y}^n_{\tau_n}, \phi \ra > \gamma) \leq \f{1}{\gamma}  c_1 \parallel \phi \parallel_\infty E(\sup_{\epsilon \leq t \leq T} \la {\mathcal Y}^n_{t}, 1\ra) \leq \f{1}{\gamma}  c_2 \parallel \phi \parallel_\infty   }

By the strong Markov Property applied to the process ${\mathcal Y}^n$
we obtain that for $\alpha, \beta \geq 0,$ we have \beas L_n(\delta_n;
\alpha, \beta) &=& E(\exp(-\alpha \la {\mathcal Y}^n_{\tau_n +
\delta_n}, \phi\ra - \beta\la {\mathcal Y}^n_{\tau_n}, \phi\ra)) \\&=&
E(\exp(-\la{\mathcal Y}^n_{\tau_n}, u^n_{\delta_n}(\alpha\phi) +
\beta\phi \ra)) \\ &=&E(\exp(-\la{\mathcal Y}^n_{\tau_n-\epsilon},
u^n_\epsilon(u^n_{\delta_n}(\alpha\phi) + \beta\phi )\ra)) \eeas
Therefore \beas \lefteqn{\mid L_n(0; \alpha, \beta) - L_n(\delta_n; \alpha,
\beta)\mid \leq } \\& \leq& \parallel
u^n_\epsilon(u^n_{\delta_n}(\alpha\phi ) + \beta\phi)-
u^n_\epsilon((\alpha + \beta)\phi) \parallel_\infty E(\sup_{t\leq T}
\la {\mathcal Y}^n_t, 1 \ra)\\ &\leq & c_1  \parallel u^n_\epsilon(u^n_{\delta_n}(\alpha\phi ) +
\beta\phi)- u^n_\epsilon((\alpha + \beta)\phi) \parallel_\infty \eeas
where is the last inequality is by Doob's maximal inequality. Now, 
\beas \lefteqn{\parallel u^n_\epsilon(u^n_{\delta_n}(\alpha\phi ) + \beta\phi)-
u^n_\epsilon((\alpha + \beta)\phi) \parallel_\infty \leq  \parallel
R^n_{n\epsilon}(u^n_{\delta_n}(\alpha\phi) - \alpha\phi)
\parallel_\infty  + }\\ && + c_2 \parallel \phi \parallel_\infty\int_0^\epsilon da \parallel u^n_{a}(
u^n_{\delta_n}(\alpha\phi) +\beta\phi)- u^n_{a}((\alpha + \beta)\phi)
\parallel_\infty  + d_n (\phi), \eeas where $d_n(\phi) = \lambda \int_0^\epsilon da \parallel \epsilon_n( u^n_{a}(
u^n_{\delta_n}(\alpha\phi) + \beta\phi) +
\epsilon_n(u^n_{a}((\alpha + \beta)\phi)) \parallel_\infty.$
Observe that \beas \lefteqn{\parallel
R^n_{n\epsilon}(u^n_{\delta_n}(\alpha\phi) - \alpha\phi)
\parallel_\infty ~~\leq~~ \parallel
R^n_{n\epsilon}(u^n_{\delta_n}(\alpha\phi) -
R^n_{n\delta_n}(\alpha\phi))\parallel_\infty } \\ && \push \push + \parallel
R^n_{n\epsilon}(R^n_{n\delta_n}(\alpha\phi)-\alpha\phi
)\parallel_\infty \\ &\leq& \parallel u^n_{\delta_n}(\alpha\phi) -
R^n_{n\delta_n}(\theta_2\phi) \parallel_\infty + \parallel
R^n_{n(\epsilon +\delta_n)}(\alpha\phi)- R^n_{n\epsilon}(\alpha\phi
)\parallel_\infty\\ &\leq & \parallel R^n_{n\delta_n}(n(1-e^{\f{s
\phi}{n}})- \alpha\phi) \parallel_\infty + \int_0^{\delta_n} da
\parallel R^n_{n(\delta_n -a)}(n^2 \Psi(\f{u^n_a \phi}{n}))
\parallel_\infty \\ && \push + \parallel R^n_{n(\epsilon
+\delta_n)}(\alpha\phi)- R^n_{n\epsilon}(\alpha\phi
)\parallel_\infty\\ &\leq& \parallel n(1-e^{\f{s \phi}{n}})-
\alpha\phi \parallel_\infty + \delta_n c_2 (\parallel \phi
\parallel^2_\infty + 1) + \parallel R^n_{n(\epsilon
+\delta_n)}(\alpha\phi)- R^n_{n\epsilon}(\alpha\phi )\parallel_\infty,\\
&& \equiv e_n(\phi) \eeas Consequently, \beas \lefteqn{\parallel
u^n_\epsilon(u^n_{\delta_n}(\alpha\phi ) + \beta\phi)-
u^n_\epsilon((r+s)\phi) \parallel_\infty \leq  e_n(\phi) + d_n (\phi)} \\
&& + c_2 \parallel \phi \parallel_\infty \int_0^\epsilon da\parallel u^n_{a}(
u^n_{\delta_n}(\alpha\phi) +\beta\phi)- u^n_{a}((r + s)\phi)
\parallel_\infty. \eeas By Proposition \ref{bcp}, $e_n(\phi) \rightarrow 0$  and $d_n(\phi) \rightarrow 0$ by our assumption
$F_n.$ Hence by a standard Gronwall argument we have that, \ben
\label{nel} \mid L_n(0; s, r) - L_n(\delta_n; s, r)\mid
\rightarrow 0 \een

By (\ref{tight1}), $\{\la {\mathcal Y}^n_{\tau_n}, \phi \ra; n =
1,2,\ldots \}$ is tight in $\R_+$. Take an arbitrary subsequence. Then
there is a further subsequence of it indexed by $\{n_k; k = 1,2,
\ldots\}$ such that $\la {\mathcal Y}^{n_k}_{\tau_{n_k}}, \phi \ra$
converges in distribution to some random limit $b$. Thus we get
$$( {\mathcal Y}^{n_k}_{\tau_{n_k}}(\phi),
{\mathcal Y}^{n_k}_{\tau_{n_k}} (\phi))\cid (b,b) \mbox{ as } k \rightarrow \infty.$$ But  (\ref{nel}) implies that
$$( {\mathcal Y}^{n_k}_{\tau_{n_k}}(\phi),
{\mathcal Y}^{n_k}_{\tau_{n_k} + \delta_{n_k}} (\phi))\cid (b,b) \mbox{ as } k \rightarrow \infty.$$

This implies that $\la {\mathcal
Y}^{n_k}_{\tau_{n_k} + \delta_{n_k}}, \phi\ra - \la {\mathcal Y}^{n_k}_{\tau_{n_k}}, \phi\ra \cid 0 \mbox{ as } k \rightarrow \infty.$  So (\ref{tight0}) holds and the proof is complete.\qed

{\bf Proof of Theorem \ref{super}} Proposition \ref{ucc} shows
that the log-Laplace functionals of the process ${\mathcal
Y}^n_t$ converge to ${\mathcal Y}_t$ for every $t \geq \epsilon
$. Proposition \ref{tight} implies tightness for the processes.
As the solution to (\ref{R}) is unique, we are done.  \qed

\end{document}